\documentclass[10pt]{article}
\usepackage{amssymb,amsthm,amsmath,wasysym}
\newtheorem{theorem}{Theorem}
\newtheorem*{thmnn}{Theorem}

\AtEndDocument{\bigskip{\footnotesize%
  \textsc{G. Durham, Department of Mathematics, University of Georgia, Athens, GA., 30602} \par  
  \textit{E-mail address}, G.~Durham: \texttt{gjdurham@uga.edu} \par
}}

\title{Representation of Integers by Ternary Quadratic Forms: A Geometric Approach}
\author{Gabriel Durham}

\begin{document}

\begin{center}
REPRESENTATION OF INTEGERS BY TERNARY QUADRATIC FORMS: A GEOMETRIC APPROACH

\textbf{\\}

GABRIEL DURHAM \footnote{This research was supported by the National Science Foundation (DMS-1461189). I would additionally like to thank Dr. Katherine Thompson,  Dr. Jeremy Rouse, Dr. Pete Clark, Sarah Blackwell and Tiffany Treece for their support. I would finally like to thank Dr. Theodore Shifrin for his guidance over the last three years and his decades of service to the University of Georgia and its students.}
\end{center}
\begin{abstract}
In 1957 N.C. Ankeny provided a new proof of the three squares theorem using geometry of numbers. This paper generalizes Ankeny's technique, proving exactly which integers are represented by $x^2 + 2y^2 + 2z^2$ and $x^2 + y^2 + 2z^2$ as well as proving sufficient conditions for an integer to be represented by $x^2+y^2+3z^2$ and $x^2 + y^2 + 7z^2$.
\end{abstract}

\section{Introduction}

A natural question in the study of quadratic forms concerns the representation of integers by certain quadratic forms. Given an $n$-ary quadratic form $Q$ and an integer $m$, does there exist a vector $\vec{x} \in\mathbb{Z}^n$ such that $Q(\vec{x}) = m$? This paper considers this question for four ternary quadratic forms: $x^2 + 2y^2 + 2z^2$, $x^2 + y^2 + 2z^2$, $x^2 + y^2 + 3z^2$, and $x^2 + y^2 + 7z^2$. We are able to prove exactly which integers are represented by the first two forms and provide sufficient conditions for an integer $m$ to be represented by the final two forms.

In 1931 Burton Jones \cite{Jones} provided one of the key breakthroughs in the study of positive definite ternary quadratic forms when he proved that the forms of a given genus collectively represent all positive integers not ruled out by certain congruence conditions. As a result, if a quadratic form is alone in its genus then one can show it represents an integer $m$ by showing that it locally represents $m$. While Jones' work made the task of determining which integers are represented by our first three forms\footnote{As well as all other ``class number one'' forms (i.e. forms which are alone in their genus).} more straightforward, the proofs presented in this paper make no use of Jones' result. Furthermore, Jones' result tells us little about forms with one or more ``genus-mates.'' While ternary quadratic forms have been studied for centuries, relatively little is known about the representation of integers by forms not alone in their genus. In particular, the question of exactly which integers are represented by $x^2+y^2+7z^2$ remains open, despite the fact that this form has the smallest determinant (and thus, in a sense, is the most simple) of any classically integral ternary quadratic form not alone in its genus. While this paper does not fully answer the question, we do provide a new proof of a result of Kaplansky \cite{Kaplansky} in our proof of Theorem $3(a)$. For a more comprehensive list of what is known about the representation of integers by $x^2 + y^2 + 7z^2$ the reader may consult Wang and Pei \cite{WP}.

To prove our results we generalize a method developed by N.C. Ankeny \cite{Ankeny}, which he used to provide a new proof of the Gauss-Legendre three squares theorem. Ankeny began with a positive integer $m$ not of the form $4^k(8\ell + 7)$ for some $k,\ell \in\mathbb{Z}$ and, using Dirichlet's theorem on primes in an arithmetic progression, defined a prime $q$ based on the prime factors of $m$. He then defined a linear map, $\vec{\Phi}:(x,y,z) \mapsto (R,S,T)$. By considering the body $\Omega = \{(R,S,T) \in\mathbb{R}^3|R^2 + S^2 + T^2 < 2m \}$ he was able to invoke Minkowski's convex body theorem to guarantee integer values of $x,y,z$ such that $\vec{\Phi}(x,y,z) \in\Omega$. By the properties of the transformation $T$ he shows that $R^2 + S^2 + T^2 = m$ and $R\in\mathbb{Z}$. To complete the proof of the three squares theorem Ankeny showed that $S^2 + T^2$ is a sum of two integer squares by showing that for all primes $p$ dividing $S^2 + T^2$ to an odd power, $\left(\frac{-1}{p}\right) = 1$. 

We alter this method by changing the transformation $\vec{\Phi}$ to show that $S^2+T^2$ is represented by other binary forms. We obtain the following results:

\begin{theorem}
The quadratic form $x^2+2y^2+2z^2$ represents a positive integer $m$ if and only if $m$ is not of the form $4^k(8\ell + 7)$ for some $k,\ell \in\mathbb{Z}$.
\end{theorem}

\begin{theorem}
The quadratic form $x^2+y^2+2z^2$ represents a positive integer $m$ if and only if $m$ is not of the form $4^k(16\ell +14)$ for some $k,\ell \in\mathbb{Z}$.
\end{theorem}

\begin{theorem}
$(a)$ If $m$ is a positive integer of the form $4^k(8\ell + 5)$ for some $k,\ell \in\mathbb{Z}$ and $\operatorname{ord}_7(m)$ is even then $m$ is represented by the quadratic form $x^2+y^2+7z^2$.

$(b)$ If $m$ is a positive integer of the form $4^k(8\ell + 1)$ for some $k,\ell \in\mathbb{Z}$ and $\operatorname{ord}_3(m)$ is even then $m$ is represented by the quadratic form $x^2+y^2+3z^2$.
\end{theorem}

\section{Background}
The following section serves as a brief introduction to the material presented in the paper.

Let $n\in\mathbb{N}$. An \underline{$n-$ary integral quadratic form}, $Q$, is a homogeneous polynomial of degree two, i.e.,
\begin{align*}
Q: \mathbb{Z}^n & \rightarrow \mathbb{Z} \\
Q: (\vec{x}) & \mapsto \sum_{1 \leq i \leq j \leq n}{a_{i,j}x_ix_j}
\end{align*}
where $a_{i,j} \in\mathbb{Z}$ for all $1 \leq i \leq j \leq n$. We say $Q$ \underline{represents} an integer $m$ if there exists an $\vec{x}\in\mathbb{Z}^n$ such that $Q(\vec{x}) = m$.

A quadratic form $Q$ is \underline{positive definite} if
\begin{align*}
(i)Q(\vec{x}) \geq 0 & \text{ for all } \vec{x}\in\mathbb{Z}^n, \\
(ii)Q(\vec{x}) = 0 & \text{ if and only if } \vec{x} = \vec{0}.
\end{align*}
Henceforth, by ``form'' we mean ``positive definite ternary integral quadratic form.''

Throughout the paper $\left(\frac{a}{b}\right)$ refers to the Jacobi symbol.

Furthermore, a quadratic form $Q$ is \underline{multiplicative} if $Q$ represents integers $a$ and $b$ implies that $Q$ represents $ab$.

The following proofs rely on two external theorems, Dirichlet's theorem on primes in an arithmetic progression and Minkowski's theorem on convex symmetric bodies. We also use the following result: For $c = 2,3$ and $7,$ $x^2 + cy^2$ represents a positive integer $m$ if and only if $\left(\frac{-c}{p}\right) = 1$ for all primes $p$ dividing $m$ to an odd power \cite{CHPT}. We end this section with the statement of Dirichlet's and Minkowski's theorems:

\begin{thmnn}(Dirichlet, 1837):
If $b,m \in\mathbb{Z}$ and $\gcd(b,m) = 1$, then there are infinitely many primes $q\equiv{b} \pmod {m}$ \cite{Apostol}.
\end{thmnn}

Before we state Minkowski's convex body theorem we will introduce the following terminology: A \underline{lattice}, $\Lambda$, in $\mathbb{R}^n$ is a subgroup of $\mathbb{R}^n$ with a basis $\{e_1, e_2, ... , e_n\}$ such that $\Lambda = \{a_1e_1 + a_2e_2 + ... + a_ne_n | a_i \in\mathbb{Z} \text{ for all } 1\leq i \leq n \}$. We call the set $T = \{b_1e_1 + b_2e_2 + ... + b_ne_n | 0 \leq b_i < 1  \text{ for all } 1\leq i \leq n \}$ the \underline{fundamental domain} of $\Lambda$ \cite{Martinet}. 
 
\begin{thmnn}(Minkowski, 1889):
Let $L$ be an $n-$dimensional lattice in $\mathbb{R}^n$ with fundamental domain $T$, and let $X$ be a bounded symmetric convex subset of $\mathbb{R}^n$. If $vol(X) > 2^nvol(T)$ then $X$ contains a non-zero point of $L$ \cite{ST}.
\end{thmnn}

\section{Proofs of Theorems}

For Theorems 1 and 2 we can look at these forms$\pmod 8$ and$\pmod{16}$ and see that they do not represent any positive integer congruent to $7 \pmod 8$ and $14 \pmod{16}$, respectively. Thus the proofs of these theorems will provide necessary and sufficient conditions for a positive integer to be represented by these forms.

We can also see that for any positive $m\in\mathbb{Z}$, if $x^2 + 2y^2 + 2z^2$ represents $m$ then $x^2 + y^2 + 2z^2$ represents $2m$, and if $x^2 + y^2 + 2z^2$ represents $m$, then $x^2 + 2y^2 + 2z^2$ represents $2m$. Thus, Theorems 1 and 2 can be proven simply by showing which positive odd integers are represented by these forms; however, the purpose of this paper is to display the manner in which Ankeny's technique can be generalized. In the spirit of this purpose we completely prove Theorem 1 using Ankeny's technique. We also use Ankeny's technique to show which positive odd integers are represented by the form $x^2 + y^2 + 2z^2 $, but, for brevity's sake, we refer to Theorem 1 when proving which even integers are represented by the form and leave the ``Ankeny-style'' proof to the reader.
\subsection{$x^2 + 2y^2 + 2z^2$}

Here we wish to prove that the quadratic form $x^2+2y^2+2z^2$ represents all positive integers not of the form $4^k(8\ell + 7)$ for some $k,\ell \in\mathbb{Z}$.

\begin{proof}
\textbf{ }

Let $m\in\mathbb{Z}$, $m>0$, where there is no $k,\ell \in\mathbb{Z}$ such that $m = 4^k(8\ell + 7)$. We want to show that $x^2+2y^2+2z^2$ represents $m$. Without loss of generality we can assume that $m$ is squarefree. We will first consider the case where $m\equiv3 \pmod 8$.

Let $q$ be an odd prime such that $q>m$, $q\equiv1 \pmod8$, and $\left(\frac{-2q}{p}\right) = 1$ for all primes $p|m$ (we know such a $q$ exists by Dirichlet's theorem on primes in an arithmetic progression). This construction of $q$ guarantees that there exists a $t\in\mathbb{Z}$ with $t^2 \equiv{ \frac{-1}{2q} } \pmod m$. This construction of $q$ also guarantees
$$
1  =  \prod_{p|m}{\left(\frac{-2q}{p}\right)}
   =  \left(\frac{-2}{m}\right)\prod_{p|m}{\left(\frac{q}{p}\right)}
   =  \prod_{p|m}{\left(\frac{q}{p}\right)}
   =  \prod_{p|m}{\left(\frac{p}{q}\right)}
   =  \left(\frac{m}{q}\right)
   =  \left(\frac{-m}{q}\right).
$$
Thus there exists a $b\in\mathbb{Z}$ (by the Chinese Remainder Theorem we can assume $b$ is odd) such that $b^2\equiv{-m} \pmod q$. As a result, there exists an $h_1 \in\mathbb{Z}$ such that $b^2 -qh_1 = -m$. Looking at this statement modulo $2$ we get $1 - qh_1 \equiv{1} \pmod{2}$ and so $-qh_1 \equiv 0 \pmod 2$. We now see that $h_1$ is even and thus $h_1 = 2h$ for some $h\in\mathbb{Z}$. Therefore we have that:
$$b^2 - 2qh = -m.$$

We now consider the body, $\Omega$, defined by:
$$\Omega = \{(R,S,T) \in\mathbb{R}^3 | 2R^2 + S^2 + T^2 < 2m\},$$
where
$$\begin{array}{ccccccc} R & = & tqx & + & bty & + & mz \\ &&&&&& \\ S & = & \ \sqrt{q}x & + & \frac{b}{\sqrt{q}}y & & \\ &&&&&& \\ T & = & & & \frac{\sqrt{m}}{\sqrt{q}}y & & \\ &&&&&&\end{array}.$$
Note that $vol(\Omega) = (\frac{4}{3})(\pi)(\frac{1}{\sqrt{2}})(\sqrt{2m}^3) = \frac{8\pi m^{3/2}}{3}.$

Let $\vec{\Phi}:\mathbb{R}^3 \rightarrow \mathbb{R}^3$, where $\vec{\Phi}:(x,y,z)\mapsto (R,S,T)$, be the function associated to the above transformation. We see that $\vec{\Phi}(\vec{x}) = M_{\vec{\Phi}}\vec{x}$ where $M_{\vec{\Phi}} = $
$$\left[ \begin{array}{ccc} tq & bt & m \\ \sqrt{q} & \frac{b}{\sqrt{q}} & 0 \\ 0 & \frac{\sqrt{m}}{\sqrt{q}} & 0 \end{array} \right].$$
By looking at $\vec{\Phi}$ as a linear map we see that $M_{\vec{\Phi}}$ is the standard matrix of $\vec{\Phi}$. Furthermore, $\det(M_{\vec{\Phi}}) = (m)(\sqrt{q})(\frac{\sqrt{m}}{\sqrt{q}}) = m^{3/2}.$

We now wish to examine $\vec{\Phi}^{-1}(\Omega)$. Since $\Omega$ is a convex symmetric body we know that $\vec{\Phi}^{-1}(\Omega)$ is as well\footnote{Since $\vec{\Phi}$ is an invertible linear map.}. We see that:
$$
vol(\vec{\Phi}^{-1}(\Omega)) = 
\frac{1}{|\det(M_{\vec{\Phi}})|} vol(\Omega) = 
\frac{1}{m^{3/2}}vol(\Omega) = 
\frac{8\pi}{3} > 8 = 2^3.
$$
By Minkowski's theorem on convex symmetric bodies there exist $x_1,y_1,z_1 \in\mathbb{Z}$, not all zero, such that $\vec{\Phi}(x_1,y_1,z_1) \in \Omega.$ Let $(R_1,S_1,T_1):=\vec{\Phi}(x_1,y_1,z_1)$.

We see that
\begin{align*}
2R_1^2 + S_1^2 + T_1^2 & \equiv 2(tqx_1 + bty_1)^2 + (\sqrt{q}x_1 + \frac{by_1}{\sqrt{q}})^2 +(\frac{\sqrt{m}y_1}{\sqrt{q}})^2\\
&\equiv 2t^2(qx_1 + by_1)^2 + \frac{1}{q}((qx_1+by_1)^2 + my_1^2)\\
&\equiv \frac{-2}{2q}(qx_1 + by_1)^2 + \frac{1}{q}(qx_1+by_1)^2\\
&\equiv 0 \pmod{m},
\end{align*}
and also that
\begin{align*}
S_1^2 + T_1^2 & = \frac{1}{q}(q^2x_1^2 + 2qbx_1y_1 + (b^2+m)y_1^2) \\
& = \frac{1}{q}(q^2x_1^2 + 2qbx_1y_1 + 2qhy_1^2) \\
& = qx_1^2 + 2bx_1y_1 + 2hy_1^2.
\end{align*}
While $S_1,T_1 \not\in\mathbb{Z}$, we can see that $S_1^2+T_1^2 \in\mathbb{Z}$. Thus $2R_1^2 + S_1^2 + T_1^2\in\mathbb{Z}$, $0<2R_1^2 + S_1^2 + T_1^2 < 2m$ and $2R_1^2 + S_1^2 + T_1^2\equiv0 \pmod{m}$. We can now conclude that $2R_1^2 + S_1^2 + T_1^2 = m$.

Our goal now is to show that $2R_1^2 + S_1^2 + T_1^2$ is represented by the quadratic form $x^2+2y^2+2z^2$. Since $R_1 \in\mathbb{Z}$ we wish to show that $S_1^2 + T_1^2$ is of the form $a^2+2b^2$ for some $a,b \in\mathbb{Z}$. It will suffice to show that $\left(\frac{-2}{p}\right)$ = 1 for all primes $p$ dividing $S_1^2 + T_1^2$ to an odd power. We note here that $x^2+2y^2$ is multiplicative and represents $2$, so we need not consider $p = 2$; furthermore, $S_1^2+T_1^2 \leq m < q$, and so we need not consider $p = q$. Thus it will suffice to show $\left(\frac{-2}{p}\right) = 1$ for all odd primes $p|v$, where $v:=q(S_1^2 + T_1^2)$, i.e. $v = (qx_1 + by_1)^2 + my_1^2$. We will do this in two cases.

\textbf{Case 1} $p\nmid m$:

We see that $0\not\equiv m \equiv 2R_1^2 \pmod p$; it follows that $-m\equiv{-2R_1^2} \pmod p$ and so $\left(\frac{-2}{p}\right) = \left(\frac{-m}{p}\right)$. Furthermore, $0 \equiv{v} \equiv{(qx_1 + by_1)^2 + my_1^2} \pmod p$ so $(qx_1 + by_1)^2 \equiv{(-m)y_1^2} \pmod p$. Since $p$ divides $v$ to an odd power we know $y_1 \not\equiv 0 \pmod p$ and thus we can conclude $1 = \left(\frac{-m}{p}\right) = \left(\frac{-2}{p}\right)$.

\textbf{}

\textbf{Case 2} $p|m$:

Since $m$ is assumed squarefree we can assume $p\nmid \frac{m}{p}$. Since $2R_1^2 + S_1^2 + T_1^2 = m$, we know $2R_1^2 \equiv{0} \pmod p$. Thus $p|R_1^2$, and so $p|R_1$. Furthermore, $(qx_1 + by_1)^2 + my^2 \equiv{(qx_1 + by_1)^2}\equiv{0} \pmod p$, so $p|(qx_1 + by_1)$.
Hence, $\frac{2R_1^2}{p} + \frac{S_1^2 + T_1^2}{p} = \frac{2R_1^2}{p} + \frac{1}{q}(\frac{(qx_1 + by_1)^2}{p} + \frac{m}{p}y_1^2) = \frac{m}{p}$. We now see that $0 + \frac{1}{q}(0 + \frac{m}{p}y_1) \equiv{\frac{m}{p}} \pmod p$. This shows that $\frac{1}{q}y_1^2 \equiv{1} \pmod p$, so $y_1^2 \equiv{q} \pmod p$ and $-2y_1^2 \equiv{-2q} \pmod p$. Thus, $\left(\frac{-2}{p}\right) = \left(\frac{-2q}{p}\right) = 1$ (by the construction of $q$).

In both cases $\left(\frac{-2}{p}\right) = 1$, thus there exist $a,b \in\mathbb{Z}$ such that $S_1^2 + T_1^2 = a^2 + 2b^2$.

\noindent Therefore, $m$ is represented by the form $x^2+2y^2+2z^2$, as required.\\

\noindent \textbf{FURTHER CASES:}

\underline{For $m$ odd}:

$\bullet$ If $m\equiv{1,5} \pmod 8$ take $q$ to be an odd prime such that $q>m$, $q\equiv{1} \pmod 8$ and $\left(\frac{-q}{p}\right) = 1$ for all primes $p|m$. We see that
$$
1  =  \prod_{p|m}{\left(\frac{-q}{p}\right)}
   =  \left(\frac{-1}{m}\right)\prod_{p|m}{\left(\frac{q}{p}\right)}
   =  \prod_{p|m}{\left(\frac{q}{p}\right)}
   =  \prod_{p|m}{\left(\frac{p}{q}\right)}
   =  \left(\frac{m}{q}\right)
   =  \left(\frac{-m}{q}\right).
$$
We now proceed as before; however, we take $t^2\equiv{\frac{-1}{4q}} \pmod m$ and take the following transformations:
$$\begin{array}{ccccccc} R & = & 2tqx & + & bty & + & mz \\ &&&&&& \\ S & = & \ \sqrt{2q}x & + & \frac{b}{\sqrt{2q}}y & & \\ &&&&&& \\ T & = & & & \frac{\sqrt{m}}{\sqrt{2q}}y & & \\ &&&&&&\end{array}.$$

\underline{For $m$ even}:

Take $m = 2m_1$. Since $m$ is assumed squarefree we can assume that $m_1$ is odd.

$\bullet$ If $m_1\equiv{1,3} \pmod 8$ take $q$ to be an odd prime such that $q\equiv{1} \pmod 8$, $q>m_1$, and $\left(\frac{-2q}{p}\right) = 1$ for all $p|m_1$. We see that
$$
1  =  \prod_{p|m_1}{\left(\frac{-2q}{p}\right)}
   =  \left(\frac{-2}{m_1}\right)\prod_{p|m}{\left(\frac{q}{p}\right)}
   =  \prod_{p|m_1}{\left(\frac{q}{p}\right)}
   =  \prod_{p|m_1}{\left(\frac{p}{q}\right)}
   =  \left(\frac{m_1}{q}\right)
   =  \left(\frac{-2m_1}{q}\right).
$$
Now take $t,b\in\mathbb{Z}$ such that $t^2 \equiv{\frac{-1}{2q}} \pmod{m_1}$ and $b^2\equiv{-2m_1} \pmod q$, where $b$ is even. Thus there exists an $h \in\mathbb{Z}$ such that $$b^2 - 2qh = -2m_1.$$ We now consider the body defined by: $$\Omega = \{(R,S,T) \in\mathbb{R}^3 | R^2 + \frac{S^2+T^2}{2} < 2m_1\}$$
(note that $vol(\Omega) = \frac{16\pi \sqrt{2} m_1^{3/2}}{3}$) where $$\begin{array}{ccccccc} R & = & tqx & + & bty & + & m_1 z \\ &&&&&& \\ S & = & \ \sqrt{q}x & + & \frac{b}{\sqrt{q}}y & & \\ &&&&&& \\ T & = & & & \frac{\sqrt{2m_1}}{\sqrt{q}}y & & \\ &&&&&&\end{array}.$$

We define $\vec{\Phi}$ to be the function associated to the above transformation. We see that $\vec{\Phi}$ is an invertible linear map and the determinant of its standard matrix is $\sqrt{2}m_1^{3/2}$, thus the volume of $\vec{\Phi}^{-1}(\Omega)$ is $\frac{16\pi}{3} > 16 = 2^3 (2)$.

We now consider the lattice $\Lambda = \{(x,y,z) \in\mathbb{Z}^3 | x\equiv{0} \pmod{2}\}$. Since \\ $[\mathbb{Z}^3:\Lambda] = 2$, we know the fundamental domain of $\Lambda$ has volume $2$ and thus we can invoke Minkowski's theorem on convex symmetric bodies to show there exists a nonzero point $(x_1,y_1,z_1) \in\Lambda$ such that $\vec{\Phi}(x_1,y_1,z_1) \in\Omega$. \\Let $(R_1,S_1,T_1):=\vec{\Phi}(x_1,y_1,z_1)$.

We see that
\begin{align*}
R_1^2 + \frac{S_1^2 + T_1^2}{2} & \equiv t^2(qx_1 + by_1)^2 + \frac{1}{2q}(qx_1 + by_1)^2 +\frac{1}{2q}({\sqrt{2m_1}y_1})^2\\
&\equiv t^2(qx_1 + by_1)^2 + \frac{1}{2q}((qx_1+by_1)^2 + 2m_1y_1^2)\\
&\equiv \frac{-1}{2q}(qx_1 + by_1)^2 + \frac{1}{2q}(qx_1+by_1)^2\\
&\equiv 0 \pmod{m_1}
\end{align*}
and that
\begin{align*}
\frac{1}{2}(S_1^2 + T_1^2) & = \frac{1}{2q}((qx_1 + by_1)^2 + 2m_1y_1^2)\\
& = \frac{1}{2q}(q^2x_1^2 + 2qbx_1y_1 + 2qhy_1^2)\\
& = q\frac{x_1^2}{2} + bx_1y_1 + hy_1^2.
\end{align*}
Since $x_1$ is even, we see that $\frac{S_1^2 + T_1^2}{2} \in\mathbb{Z}$. Since $R_1 \in\mathbb{Z}$ as well, $R_1^2 + \frac{S_1^2 + T_1^2}{2} \in\mathbb{Z}$. Moreover, $0< R_1^2 + \frac{S_1^2 + T_1^2}{2} < 2m_1$, thus we conclude that $R_1^2 + \frac{S_1^2 + T_1^2}{2} = m_1$.

We now define $v:=2q(S_1^2 + T_1^2)$, i.e. $v = (qx_1+by_1)^2 + 2m_1y_1^2$. We now wish to show that $\frac{S_1^2+T_1^2}{2}$ is represented by the binary quadratic form $x^2+2y^2$. As before, it will suffice to show $\left(\frac{-2}{p}\right) = 1$ for all odd primes $p$ dividing $v$ to an odd power. We proceed in two cases.

\textbf{Case 1} $p\nmid m$:

Since $m_1\equiv{R_1^2} \pmod p$, $\left(\frac{m_1}{p}\right) = 1$. Furthermore, $(qx_1 +by_1)^2 \equiv -2m_1y_1^2 \pmod p$. It follows that $\left(\frac{-2}{p}\right) = \left(\frac{m_1}{p}\right) = 1$.

\textbf{Case 2} $p|m$:

By the same argument presented in the $m\equiv3 \pmod 8$ case we see that $\frac{1}{2q}(2y_1^2)\equiv1 \pmod p$ and so $-2y_1^2\equiv{-2q} \pmod p$. By the construction of $q$, $\left(\frac{-2}{p}\right) = \left(\frac{-2q}{p}\right) = 1$.

In both cases we have $\left(\frac{-2}{p}\right) = 1$. Therefore there exist $a_1,b_1 \in\mathbb{Z}$ such that $\frac{S_1^2+T_1^2}{2} = a_1^2 + 2b_1^2$. Since the form $x^2+2y^2$ is multiplicative and represents $2$ we also know there exist $a,b\in\mathbb{Z}$ such that $S_1^2 + T_1^2 = a^2 + 2b^2$. Furthermore, since $R_1^2 + \frac{S_1^2 + T_1^2}{2} = m_1$, $2R_1^2 + S_1^2 + T_1^2 = m$. Noting that $R_1\in\mathbb{Z}$ we see $m = 2R_1^2 + a^2 + 2b^2$.

Therefore $m$ is represented by the quadratic form $x^2+2y^2+2z^2$, as required.

$\bullet$ If $m_1\equiv5 \pmod 8$, take $q$ to be an odd prime such that $q>m_1$, $q\equiv 5 \pmod 8$, and $\left(\frac{-2q}{p}\right) = 1$ for all primes $p|m_1$. We see that 
$$
1  =  \prod_{p|m_1}{\left(\frac{-2q}{p}\right)}
   =  \left(\frac{-2}{m_1}\right)\prod_{p|m}{\left(\frac{q}{p}\right)}
   =  (-1) \prod_{p|m_1}{\left(\frac{q}{p}\right)}
   =  (-1) \prod_{p|m_1}{\left(\frac{p}{q}\right)}
   =  (-1) \left(\frac{m_1}{q}\right)
   =  \left(\frac{-2m_1}{q}\right).
$$
The rest of the proof proceeds as above.

$\bullet$ If $m_1\equiv7 \pmod 8$, take $q$ to be an odd prime such that $q>m_1$, $q\equiv 3 \pmod 8$, and $\left(\frac{-2q}{p}\right) = 1$ for all primes $p|m_1$. We see that 
$$
1  =  \prod_{p|m_1}{\left(\frac{-2q}{p}\right)}
   =  \left(\frac{-2}{m_1}\right)\prod_{p|m}{\left(\frac{q}{p}\right)}
   =  (-1) \prod_{p|m_1}{\left(\frac{q}{p}\right)}
   =  \prod_{p|m_1}{\left(\frac{p}{q}\right)}
   =  \left(\frac{m_1}{q}\right)
   =  \left(\frac{-2m_1}{q}\right).
$$
The rest of the proof proceeds as above.

This completes the proof of Theorem 1.
\end{proof}

\subsection{$x^2 + y^2 + 2z^2$}

We now wish to prove that the quadratic form $x^2+y^2+2z^2$ represents all positive integers not of the form $4^k(16\ell +14)$ for some $k,\ell \in\mathbb{Z}$

\begin{proof}
\textbf{\\}
Let $m\in\mathbb{Z},$ where $m>0$ and there is no $k,\ell \in\mathbb{Z}$ such that $m = 4^k(16\ell +14)$. We wish to prove that $m$ is represented by the form $x^2+y^2+2z^2$. Without loss of generality we can assume that $m$ is squarefree. As before we first consider the case where $m\equiv 3 \pmod 8$.

Let $q$ be an odd prime where $q>m$, $q\equiv1 \pmod 8$ and $\left(\frac{-2q}{p}\right) = 1$ for all primes $p|m$ (Dirichlet's theorem regarding primes in an arithmetic progression guarantees us such a $q$). Furthermore, this construction of $q$ ensures that there exists a $t\in\mathbb{Z}$ such that $t^2\equiv{\frac{-1}{2q}} \pmod m$. We also have:
$$
1  =  \prod_{p|m}{\left(\frac{-2q}{p}\right)}
   =  \left(\frac{-2}{m}\right)\prod_{p|m}{\left(\frac{q}{p}\right)}
   =  \prod_{p|m}{\left(\frac{q}{p}\right)}
   =  \prod_{p|m}{\left(\frac{p}{q}\right)}
   =  \left(\frac{m}{q}\right)
   =  \left(\frac{-2m}{q}\right).
$$
Thus there exists a $b\in\mathbb{Z}$ (where $b$ is even) such that $b^2\equiv{-2m} \pmod q$. As a result we know $b^2-qh_1 = -2m$ for some $h_1\in\mathbb{Z}$. Since $b$ is even $h_1$ must be as well. We write $h_1 = 2h$ for some $h\in\mathbb{Z}$ and see $$b^2-2qh = -2m.$$

We now consider the body defined by:
$$\Omega = \{(R,S,T)\in{\mathbb{R}^3}|R^2+S^2+T^2 < 2m\}$$
(note that $vol(\Omega) = \frac{4}{3}\pi \sqrt{2m}^3 = \frac{8\sqrt{2} \pi m^{3/2}}{3}$) where:
$$\begin{array}{ccccccc} R & = & 2tqx & + & bty & + & mz \\ &&&&&& \\ S & = & \ \sqrt{2q}x & + & \frac{b}{\sqrt{2q}}y & & \\ &&&&&& \\ T & = & & & \frac{\sqrt{2m}}{\sqrt{2q}}y & & \\ &&&&&&\end{array}.$$

Let $\vec{\Phi}:\mathbb{R}^3 \rightarrow \mathbb{R}^3$, where $\vec{\Phi}: (x,y,z) \mapsto (R,S,T)$ be the function associated to the above transformation. Thus, $\vec{\Phi}(\vec{x}) = M_{\vec{\Phi}}\vec{x}$ where $M_{\vec{\Phi}} = $
$$\left[ \begin{array}{ccc} 2tq & bt & m \\ \sqrt{2q} & \frac{b}{\sqrt{2q}} & 0 \\ 0 & \frac{\sqrt{2m}}{\sqrt{2q}} & 0 \end{array} \right].$$
Since $\Omega$ is a convex symmetric body, $\vec{\Phi}^{-1}$ is as well. We see that
$$vol(\vec{\Phi}^{-1}(\Omega)) = \frac{vol(\Omega)}{|\det(M_{\vec{\Phi}})|} = \frac{8\pi}{3} > 8.$$
Minkowski's theorem on convex symmetric bodies guarantees the existence of $x_1,y_1,z_1 \in\mathbb{Z}$ (not all zero) such that $\vec{\Phi}(x_1,y_1,z_1) \in\Omega$. Let $(R_1,S_1,T_1):=\vec{\Phi}(x_1,y_1,z_1)$.
We see that,
\begin{align*}
R_1^2 + S_1^2 + T_1^2 & \equiv (2tqx_1 + bty_1)^2 + (\sqrt{2q}x_1 + \frac{by_1}{\sqrt{2q}})^2 +(\frac{\sqrt{2m}y_1}{\sqrt{2q}})^2\\
&\equiv t^2(2qx_1 + by_1)^2 + \frac{1}{2q}((2qx_1+by_1)^2 + 2my_1^2)\\
&\equiv t^2(2qx_1 + by_1)^2 + \frac{1}{2q}(2qx_1+by_1)^2\\
&\equiv 0 \pmod{m},
\end{align*}
and that
\begin{align*}
S_1^2 + T_1^2 & = \frac{1}{2q}(4q^2x_1^2 + 4qbx_1y_1 + (b^2+2m)y_1^2) \\
& = \frac{1}{2q}(4q^2x_1^2 + 4qbx_1y_1 + 2qhy_1^2) \\
& = 2qx_1^2 + 2bx_1y_1 + hy_1^2.
\end{align*}
We now see that $S_1^2 + T_1^2 \in\mathbb{Z}$ and, as before, we can conclude that $R_1^2 + S_1^2 + T_1^2 = m$.

Let $v:=2q(S_1^2+T_1^2)$, i.e. $v = (2qx_1+by_1)^2 + 2my_1^2$. We wish to show that $S_1^2 + T_1^2$ is of the form $a^2+2b^2$ for some $a,b \in\mathbb{Z}$. As was the case in our proof of Theorem 1, it will suffice to show $\left(\frac{-2}{p}\right) = 1$ for all odd primes $p$ dividing $v$ to an odd power (and we need not consider $p=q$). We will do so in two cases.

\textbf{Case 1} $p \nmid m$:

Since $m\equiv{R_1}^2 \pmod p$, $\left(\frac{m}{p}\right) = 1$. Furthermore, $(2qx_1+by_1)^2 \equiv{-2my_1^2} \pmod p$. Thus $\left(\frac{-2}{p}\right) = \left(\frac{m}{p}\right) = 1$.

\textbf{Case 2} $p|m$:

Since $m$ is assumed squarefree we know that $p\nmid \frac{m}{p}$. Since $R_1^2 + S_1^2 + T_1^2 = m$, $R_1^2 \equiv 0 \pmod p$ and so $p|R_1$. Moreover, we know that \\$(2qx_1+by_1)^2 -2my_1^2 \equiv{(2qx_1+by_1)^2}\equiv{0} \pmod p$. Thus $p|(2qx_1+by_1)$. Noting that $\frac{R_1^2}{p} + \frac{1}{2q}(\frac{(2qx_1+by_1)^2}{p} + 2\frac{m}{p}y_1^2) = \frac{m}{p}$, we see that $\frac{1}{2q}(2\frac{m}{p}y_1^2) \equiv{\frac{m}{p}} \pmod p$, and so $-2y_1^2 \equiv{-2q} \pmod p$. Thus $\left(\frac{-2}{p}\right) = \left(\frac{-2q}{p}\right) = 1$.

In both cases we have $\left(\frac{-2}{p}\right) = 1$ for all odd primes $p$ dividing $v$ to an odd power. We now know there exist $a,b\in\mathbb{Z}$ such that $S_1^2 + T_1^2 = a + 2b^2$.

\noindent Therefore, $m$ is represented by the quadratic form $x^2 + y^2 + 2z^2$, as required.

\noindent \textbf{FURTHER CASES:}

For $m$ odd:

$\bullet$ For $m\equiv 7 \pmod 8$ take $q$ to be an odd prime such that $q>m$, $q\equiv3 \pmod 8$ and $\left(\frac{-2q}{p}\right) = 1$ for all primes $p|m$. We see that
$$
1  =  \prod_{p|m}{\left(\frac{-2q}{p}\right)}
   =  \left(\frac{-2}{m}\right)\prod_{p|m}{\left(\frac{q}{p}\right)}
   =  (-1) \prod_{p|m}{\left(\frac{q}{p}\right)}
   =  \prod_{p|m}{\left(\frac{p}{q}\right)}
   =  \left(\frac{m}{q}\right)
   =  \left(\frac{-2m}{q}\right).
$$
The rest of the proof proceeds as above.

$\bullet$ For $m\equiv{1,5} \pmod 8$ take $q$ to be an odd prime such that $q>m$, $q\equiv 1 \pmod 8$ and $\left(\frac{-q}{p}\right) = 1$ for all primes $p|m$. We see that
$$
1  =  \prod_{p|m}{\left(\frac{-q}{p}\right)}
   =  \left(\frac{-1}{m}\right)\prod_{p|m}{\left(\frac{q}{p}\right)}
   =  \prod_{p|m}{\left(\frac{q}{p}\right)}
   =  \prod_{p|m}{\left(\frac{p}{q}\right)}
   =  \left(\frac{m}{q}\right)
   =  \left(\frac{-2m}{q}\right).
$$
We let $t^2 \equiv{\frac{-1}{q}} \pmod{m}$ and take the following transformation:
$$\begin{array}{ccccccc} R & = & tqx & + & bty & + & mz \\ &&&&&& \\ S & = & \ \sqrt{q}x & + & \frac{b}{\sqrt{q}}y & & \\ &&&&&& \\ T & = & & & \frac{\sqrt{2m}}{\sqrt{q}}y & & \\ &&&&&&\end{array}.$$
The rest of the proof proceeds as above.

For $m$ even:

Take $m = 2m_1$. By Theorem 1 we know, for $m_1 \equiv{1,3,5} \pmod 8$, there exist $x,y,z \in\mathbb{Z}$ such that $m_1 = x^2 + 2y^2 + 2z^2$. Thus $2m_1 = 2x^2 + (2y)^2 + (2z)^2$. Therefore, $m$ is represented by the quadratic form $x^2+y^2+2z^2$.

This completes the proof of Theorem 2.
\end{proof}

\subsection{$x^2 + y^2 + 7z^2$ and $x^2 + y^2 + 3z^2$}

Here we wish to show that if $m$ is a positive integer of the form $4^k(8\ell + 5)$ for some $k,\ell \in\mathbb{Z}$ and $\operatorname{ord}_7(m)$ is even then $m$ is represented by the quadratic form $x^2+y^2+7z^2$ and if $m$ is a positive integer of the form $4^k(8\ell + 1)$ for some $k,\ell \in\mathbb{Z}$ and $\operatorname{ord}_3(m)$ is even then $m$ is represented by the quadratic form $x^2+y^2+3z^2$.

\begin{proof}
\textbf{  }

We first seek to prove Theorem 3(a). Without loss of generality we can assume $m$ squarefree (and consequently that $7\nmid m$) and $m\equiv 5 \pmod 8$.

Let $q$ be an odd prime such that $q>m$, $q\equiv 1 \pmod{28}$ and $\left(\frac{-q}{p}\right) = 1$ for all primes $p|m$ (Dirichlet's theorem of primes in an arithmetic progression guarantees the existence of $q$). By this construction of $q$, there exists a $t\in\mathbb{Z}$ such that $t^2\equiv{\frac{-1}{4q}} \pmod m$. Additionally we see that
$$
1  =  \prod_{p|m}{\left(\frac{-q}{p}\right)}
   =  \left(\frac{-1}{m}\right)\prod_{p|m}{\left(\frac{q}{p}\right)}
   =  \prod_{p|m}{\left(\frac{q}{p}\right)}
   =  \prod_{p|m}{\left(\frac{p}{q}\right)}
   =  \left(\frac{m}{q}\right)
   =  \left(\frac{-7m}{q}\right).
$$
Thus there is a $b\in\mathbb{Z}$ (where $b$ is odd) such that $b^2 \equiv{-7m} \pmod q$. This shows us that $b^2 - qh_1 = -7m$ for some $h_1 \in\mathbb{Z}$. Looking at this statement modulo $4$ we see that $b^2 - qh_1 \equiv m \pmod 4$ and so $h_1 \equiv{0} \pmod 4$. Therefore $h_1 = 4h$ for some $h\in\mathbb{Z}$ and $$b^2-4qh = -7m.$$ Looking at the above statement modulo $8$ we see $4qh \equiv 4 \pmod 8$, thus $h$ is odd.

We now consider the body defined by:
$$\Omega = \{(R,S,T)\in{\mathbb{R}^3}|R^2+S^2+T^2 < 2m\}$$
(note that $vol(\Omega) = \frac{4}{3}\pi \sqrt{2m}^3 = \frac{8\sqrt{2} \pi m^{3/2}}{3}$) where:
$$\begin{array}{ccccccc} R & = & 2tqx & + & bty & + & mz \\ &&&&&& \\ S & = & \ \sqrt{q}x & + & \frac{b}{2\sqrt{q}}y & & \\ &&&&&& \\ T & = & & & \frac{\sqrt{7m}}{2\sqrt{q}}y & & \\ &&&&&&\end{array}.$$

Let $\vec{\Phi}:\mathbb{R}^3 \rightarrow \mathbb{R}^3$, where $\vec{\Phi}: (x,y,z) \mapsto (R,S,T)$ be the function associated to the above transformation. Thus, $\vec{\Phi}(\vec{x}) = M_{\vec{\Phi}}\vec{x}$ where $M_{\vec{\Phi}} = $
$$\left[ \begin{array}{ccc} 2tq & bt & m \\ \sqrt{q} & \frac{b}{2\sqrt{2q}} & 0 \\ 0 & \frac{\sqrt{7m}}{2\sqrt{q}} & 0 \end{array} \right].$$

We see that $\det(M_{\vec{\Phi}}) = \frac{\sqrt{7}m^{3/2}}{2}$. Since $\Omega$ is a convex symmetric body, $\vec{\Phi}^{-1}(\Omega)$ is as well. It follows that:
$$vol(\vec{\Phi}^{-1}(\Omega)) = \frac{vol(\Omega)}{|\det(M_{\vec{\Phi}})|} = \frac{16\sqrt{2} \pi}{3\sqrt{7}} > 8.$$
Minkowski's theorem on convex symmetric bodies guarantees that there exist $x_1,y_1,z_1 \in\mathbb{Z}$,  not all zero, such that $\vec{\Phi}(x_1,y_1,z_1) \in\Omega$. Let $(R_1,S_1,T_1):=\vec{\Phi}(x_1,y_1,z_1)$.
We see that
\begin{align*}
R_1^2 + S_1^2 + T_1^2 & \equiv (2tqx_1 + bty_1)^2 + (\sqrt{q}x_1 + \frac{by_1}{2\sqrt{q}})^2 +(\frac{\sqrt{7m}y_1}{2\sqrt{q}})^2\\
&\equiv t^2(2qx_1 + by_1)^2 + \frac{1}{4q}((2qx_1+by_1)^2 + 7my_1^2)\\
&\equiv t^2(2qx_1 + by_1)^2 + \frac{1}{4q}(2qx_1+by_1)^2\\
&\equiv 0 \pmod{m},
\end{align*}
and that
\begin{align*}
S_1^2 + T_1^2 & = \frac{1}{4q}(4q^2x_1^2 + 4qbx_1y_1 + (b^2+7m)y_1^2) \\
& = \frac{1}{4q}(4q^2x_1^2 + 4qbx_1y_1 + 4qhy_1^2) \\
& = qx_1^2 + bx_1y_1 + hy_1^2.
\end{align*}
Thus $S_1^2 + T_1^2 \in\mathbb{Z}$ and, as before, we conclude $R_1^2 + S_1^2 + T_1^2 = m$. Since $R_1\in\mathbb{Z}$ we wish to show $S_1^2 + T_1^2$ is of the form $a^2 + 7b^2$ for some $a,b \in\mathbb{Z}$. It will be sufficient to show that $\left(\frac{-7}{p}\right) = 1$ for all primes $p$ dividing $S_1^2 + T_1^2$ to an odd power.

Note that $p\leq m < q$, so we can assume $p\neq q$. Furthermore, if $2|S_1^2 + T_1^2 = qx_1^2 + bx_1y_1 + hy_1^2$ then $2|x,y$ (since $h$ and $b$ are odd), thus $2$ divides $S_1^2 + T_1^2$ to an even power and so we need not consider $p = 2$. We also need not consider $p = 7,$ as the form $x^2+7y^2$ represents $7$ and is multiplicative. As a result it will suffice to show $\left(\frac{-7}{p}\right) = 1$ for all primes $p\neq 2,7,q$ dividing $v$, where $v:=4q(S_1^2+T_1^2)$, i.e. $v = (2qx_1+by_1)^2 + 7my_1^2$. We proceed in two cases.

\textbf{Case 1} $p\nmid m$:

We see that $(2qx_1+by_1)^2\equiv{-7my_1^2} \pmod p$ and $m\equiv{R_1}^2 \pmod p$. It follows that $\left(\frac{-7}{p}\right) = \left(\frac{m}{p}\right) = 1$.

\textbf{Case 2} $p|m$:

Note that since $m$ is assumed squarefree $p\nmid \frac{m}{p}$. Thus, $R_1^2\equiv{0} \pmod p$ and so $p|R_1$. Additionally, $(2qx_1+by_1)^2 + 7my_1^2\equiv{(2qx_1+by_1)^2}\equiv{0} \pmod p$ so $p|(2qx_1+by_1)$. We now see that $\frac{R_1^2}{p} + \frac{1}{4q}(\frac{(2qx_1+by_1)^2}{p} + 7\frac{m}{p}y_1^2) = \frac{m}{p}$. Thus $\frac{1}{4q}(7\frac{m}{p}y_1^2)\equiv{\frac{m}{p}} \pmod p$. It follows that $-7y_1^2 \equiv{-4q} \pmod p$. By the construction of $q$, $\left(\frac{-7}{p}\right) = \left(\frac{-q}{p}\right) = 1$.

In both cases we see $\left(\frac{-7}{p}\right) = 1$. Thus there exist $a,b\in\mathbb{Z}$ such that $m = a^2 + 7b^2$.

\noindent Therefore, $m$ is represented by the quadratic form $x^2 + y^2 + 7z^2$, as required.

This completes the proof of Theorem 3(a).

\textbf{ }

The proof of Theorem 3(b) follows the above proof; however, we take $q$ to be an odd prime such that $q>m$, $q\equiv 1 \pmod{12}$ and $\left(\frac{-q}{p}\right) = 1$ for all primes $p|m$. It follows that
$$
1  =  \prod_{p|m}{\left(\frac{-q}{p}\right)}
   =  \left(\frac{-1}{m}\right)\prod_{p|m}{\left(\frac{q}{p}\right)}
   =  \prod_{p|m}{\left(\frac{q}{p}\right)}
   =  \prod_{p|m}{\left(\frac{p}{q}\right)}
   =  \left(\frac{m}{q}\right)
   =  \left(\frac{-3m}{q}\right)
$$
and we take the following transformation:
$$\begin{array}{ccccccc} R & = & 2tqx & + & bty & + & mz \\ &&&&&& \\ S & = & \ \sqrt{q}x & + & \frac{b}{2\sqrt{q}}y & & \\ &&&&&& \\ T & = & & & \frac{\sqrt{3m}}{2\sqrt{q}}y & & \\ &&&&&&\end{array}.$$

This completes the proof of Theorem 3.

\end{proof}

\end{document}